# A Noncommutative Version of the Natural Numbers

Tyler Foster

**Abstract.** — In this note, we construct and study an algebraic system similar to the natural numbers, but with noncommutative addition. The addition we introduce is a binary operation that commutes with itself in the sense of N. Durov. Neverheless, the multiplication in this system (defined by iterating the noncommutative addition) turns out to be associative, commutative, and distributive over addition, and the resulting system has interesting and nontrivial arithmetic.

**1 — Introduction.**

The operation of addition in the natural numbers $\mathbb{N} = \{1, 2, 3, \ldots\}$ is of course commutative and associative. We will construct a "noncommutative $\mathbb{N}$" by weakening this pair of properties, replacing it with a single property called *commutativity of a binary operation with itself*. It is an instance of a more general property that N. Durov introduces in [1], called *commutativity* of $m$- and $n$-ary operations:

**1.1 Definition.** — If $X$ is a set equipped with $m$- and $n$-ary operations $\varphi : X^m \to X$ and $\psi : X^n \to X$, respectively, we say that $\varphi$ and $\psi$ *commute* when

$$\varphi(\psi(x_{11}, \ldots, x_{1n}), \ldots, \psi(x_{m1}, \ldots, x_{mn})) = \psi(\varphi(x_{11}, \ldots, x_{m1}), \ldots, \varphi(x_{1n}, \ldots, x_{mn}))$$

for any choice of $x_{ij} \in X$.

**1.2** If we write the general element of $X^{mn}$ as an $m \times n$-array

$$\begin{vmatrix} x_{11} & x_{12} & \cdots & x_{1n} \\ x_{21} & x_{22} & \cdots & x_{2n} \\ \vdots & \vdots & \ddots & \vdots \\ x_{m1} & x_{m2} & \cdots & x_{mn} \end{vmatrix}$$

and think of $\varphi$ as taking columns as arguments while $\psi$ takes rows as arguments, then commutativity of $\varphi$ with $\psi$ says that we get the same value in $X$ if we **(i)** first apply $\psi$ to each row and then apply $\varphi$ to the resulting column, or **(ii)** first apply $\varphi$ to each column and then apply $\psi$ to the resulting row.

**1.3 Example.** — Our prototypical example of commuting operations is that of finite linear combinations over an associative, commutative ring $R$, that is, operations

$$\varphi(x_1, \ldots, x_m) = a_1 x_1 + \cdots + a_m x_m$$
$$\psi(y_1, \ldots, y_n) = b_1 y_1 + \cdots + b_n y_n$$

where $a_1, \ldots, a_m, b_1, \ldots, b_n \in R$, and where our set $X$ is any $R$-module. For this reason, we're to think of systems of commuting operations as generalizing systems of finite linear combinations. That is, such systems are supposed to generalize addition over a ring.

**1.4** In the present study, we focus on <u>binary</u> operations. Specifically, we focus on a single binary operation

$$\oplus : X^2 \longrightarrow X$$

that commutes with itself. Writing the definition out explicitly in this case, we see that the condition that $\oplus$ commute with itself is nothing but the condition that

$$(w \oplus x) \oplus (y \oplus z) = (w \oplus y) \oplus (x \oplus z). \qquad (*)$$

Said differently, $\oplus$ commutes with itself as long as we can transpose $x$ and $y$ in any expression $(w \oplus x) \oplus (y \oplus z)$.

We will call this property *self-commutativity of $\oplus$*.

## 2 — The Construction of "Noncommutative $\mathbb{N}$."

We want to consider the free algebraic system $A$ generated by a single symbol "1" and a single, self-commuting binary operation $\oplus$. If a self-commuting, binary operation is a generalization of addition, then a tentative analogy between this system and $\mathbb{N}$ is already clear. We will see that the analogy is even stronger than one might expect.

**2.1** We call the expression $x \oplus y$ the (*noncommutative*) *sum* of $x$ and $y$, and we let $M$ be the set of all nonassociative, bracketed sums of a single element 1 (what Bourbaki would call the *free magma* on a single generator 1). We can identify $M$ with the set of binary rooted trees, identifying each of its subsets

$$M_n := \{\text{sums containing } n\text{-instances of } 1\}$$

with the set of binary rooted trees with $n$-branches. We point out that $|M_n|$ = the $n^{\text{th}}$ Catalan number $C(n)$.

**2.2** We now obtain the system $A$ we're after by imposing the relation of self-commutativity on $M$. Specifically, we define

$$A := M / \sim,$$

where $\sim$ is the minimum equivalence relation compatible with $\oplus$ in which the identity ($*$) of §1 above holds for all $w, x, y, z \in M$. We refer to the identifications in $\sim$ given by ($*$) as *elementary instances of self-commutativity*, and all secondary identifications in $\sim$ as *algebraic consequences* (*of self-commutativity*).

**2.3** Our system $A$ consists of elements like

$$1, \quad 1 \oplus 1, \quad 1 \oplus (1 \oplus 1), \quad (1 \oplus (1 \oplus 1)) \oplus (1 \oplus 1),$$

with various identities holding between such expressions. For instance, we have the equality

$$(1 \oplus (1 \oplus 1)) \oplus (1 \oplus 1)$$
$$= (1 \oplus 1) \oplus ((1 \oplus 1) \oplus 1)$$

in $A$. Here the pair of arrows indicates an elementary instance of self-commutativity.

**2.4** If we let the *magnitude n* of an expression $a$ in $M$ be the number of 1's appearing in it, then we can make a small table listing all the elements in $A$ coming from expressions of magnitude $\leq 5$:

| 1 | 2 | 3 | 4 | 5 | |
|---|---|---|---|---|---|
| 1 | $2 := 1 \oplus 1$ | $3_1 := 1 \oplus (1 \oplus 1) = 1 \oplus 2$<br>$3_2 := (1 \oplus 1) \oplus 1 = 2 \oplus 1$ | $4_1 := 2 \oplus 2$<br>$4_2 := 1 \oplus 3_1$<br>$4_3 := 3_1 \oplus 1$<br>$4_4 := 1 \oplus 3_2$<br>$4_5 := 3_2 \oplus 1$ | $5_1 := 2 \oplus 3_1$<br>$5_2 := 3_1 \oplus 2 = 2 \oplus 3_2$<br>$5_3 := 3_2 \oplus 2$<br>$5_4 := 1 \oplus 4_2$<br>$5_5 := 4_2 \oplus 1$<br>$5_6 := 1 \oplus 4_3$<br>$5_7 := 4_3 \oplus 1$<br>$5_8 := 1 \oplus 4_4$<br>$5_9 := 4_4 \oplus 1$<br>$5_{10} := 1 \oplus 4_5$<br>$5_{11} := 4_5 \oplus 1$<br>$5_{12} := 1 \oplus 4_1$<br>$5_{13} := 4_1 \oplus 1$ | **(table 1)** |

Here, in order to avoid long, hard-to-look-at expressions in 1's and $\oplus$'s, we've named the elements recursively. Thus the identity $3_1 \oplus 2 = 2 \oplus 3_2$ appearing in the table is nothing but the identity $(1 \oplus (1 \oplus 1)) \oplus (1 \oplus 1) = (1 \oplus 1) \oplus ((1 \oplus 1) \oplus 1)$ in $A$.

As the magnitude increases, the number of identities arising between expressions for elements in $A$ increases. For small $n$, the number $D(n)$ of elements in $A$ coming from expressions of magnitude $n$ matches the Catalan number $C(n)$, but begins to diverge from $C(n)$ once $n \geq 5$. The table below exhibits this phenomena for $1 \leq n \leq 7$:

| magnitude $n$: | 1 | 2 | 3 | 4 | 5 | 6 | 7 |
|---|---|---|---|---|---|---|---|
| Catalan numbers $C(n)$: | 1 | 1 | 2 | 5 | 14 | 42 | 132 |
| # of distinct images in $A$: | 1 | 1 | 2 | 5 | 13 | 36 | 102 |

(table 2)

As an example, among all expressions of magnitude 6 in $M$ we have identities

$$2 \oplus 4_3 = 4_2 \oplus 2, \quad 2 \oplus 4_1 = 3_1 \oplus 3_1, \quad 1 \oplus (3_1 \oplus 2) = 1 \oplus (2 \oplus 3_2),$$
$$2 \oplus 4_5 = 4_4 \oplus 2, \quad 4_1 \oplus 2 = 3_2 \oplus 3_2, \quad \text{and} \quad (3_1 \oplus 2) \oplus 1 = (2 \oplus 3_2) \oplus 1,$$

accounting for the drop from the Catalan number $C(6) = 42$ to $D(6) = 36$ in the magnitude 6 column above.

I do not know a closed formula for the $D(n)$.

## 3 — Algebra in $A$.

By construction, the binary operation $\oplus$ on $M$ induces a well defined binary operation $\oplus : A \times A \to A$ that commutes with itself but is neither associative nor commutative.

**3.1** The system $M$ also comes with a natural operation of multiplication. Using our identification of $M$ with the set of binary rooted trees, we define the product $a \cdot b$ of elements $a, b \in M$ to be the expression in $M$ associated to the tree gotten by grafting a copy of $b$ onto each branch of $a$.

For instance, when $a$ is the expression $1 \oplus 1$ and $b$ is the expression $1 \oplus (1 \oplus 1)$, the grafting describing the product $a \cdot b$ is so $a \cdot b = (1 \oplus (1 \oplus 1)) \oplus (1 \oplus (1 \oplus 1))$. This generalizes the elementary school definition of multiplication in $\mathbb{N}$, namely that $m \cdot n$ is the $m$-term sum $n + \cdots + n$.

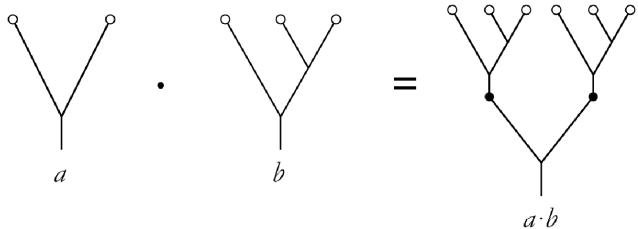

It is clear that our multiplication in $M$ at least distributes over $\oplus$ from the right.

**3.2 Claim.** — Multiplication in $M$ descends to an associative operation in $A$ with unit 1.
**Proof:** That multiplication in $M$ descends to a well defined operation in $A$ follows from the definitions of multiplication in $M$ and the equivalence relation defining $A$. Associativity is also clear, it being nothing but associativity of composition of functions, and our definition of multiplication makes 1 the obvious unit. ∎

**3.3** Observe that the identity $(1 \oplus (1 \oplus 1)) \oplus (1 \oplus (1 \oplus 1)) = (1 \oplus 1) \oplus ((1 \oplus 1) \oplus (1 \oplus 1))$ in $A$ implies that $a \cdot b = b \cdot a$ in the example of 3.1 above. In fact, one consequence of self-commutativity of $\oplus$ and our definition of multiplication is that all products commute in $A$, and multiplication distributes over $\oplus$ from the left. These facts follow from the simple principle that operations built out of commuting operations commute, which we make explicit in the following Lemma.

**3.4 Lemma.** — Let $X$ be a set and let $\Phi$ be any collection of operations $\varphi_i : X^{n_i} \to X$ (of varying arity $n_i$, respectively), all of which commute with one another in the sense of Definition 1.1. If

$$\psi_a : X^{m_1} \to X \quad \text{and} \quad \psi_b : X^{m_2} \to X$$

are each arbitrary compositions of operations $\varphi_i \in \Phi$, then $\psi_a$ and $\psi_b$ commute.
**Proof:** The claim is trivial if $\psi_a, \psi_b \in \Phi$. We proceed first by induction on the number of factors $\varphi_i \in \Phi$ in the composition forming $\psi_a$. Suppose every composition built out of $n$-many factors commutes with every $\varphi \in \Phi$, and let $\psi_a$ be built out of $(n+1)$-many factors Then $\psi_a$ is of the form

$$\psi_a(t_{11}, \ldots, t_{n_1 1}, \ldots, t_{1\,m}, \ldots, t_{n_m m}) = \varphi_i\big(\psi_1(t_{11}, \ldots, t_{n_1 1}), \ldots, \psi_m(t_{1\,m}, \ldots, t_{n_m m})\big).$$

Fix a $k$-ary operation $\varphi \in \Phi$, let $N := n_1 + \cdots + n_{m-1}$, and let $M := n_1 + \cdots + n_{m-1} + n_m$. Then for any $k \times M$-array

$$\begin{vmatrix} x_{11} & \cdots & x_{n_1 1} & \cdots & & \cdots & x_{N1} & \cdots & x_{M1} \\ \vdots & \ddots & \vdots & & & & \vdots & \ddots & \vdots \\ x_{1k} & \cdots & x_{n_1 k} & \cdots & & \cdots & x_{N1} & \cdots & x_{M1} \end{vmatrix}$$

we have $\psi_a\big(\varphi(x_{11}, \ldots, x_{1k}), \ldots, \varphi(x_{M1}, \ldots, x_{Mk})\big) =$

$$= \varphi_i\big(\psi_1(\varphi(x_{11}, \ldots, x_{1k}), \ldots, \varphi(x_{n_1 1}, \ldots, x_{n_1 k})), \ldots, \psi_m(\varphi(x_{N1}, \ldots, x_{Nk}), \ldots, \varphi(x_{M1}, \ldots, x_{Mk}))\big)$$

$$= \varphi_i\big(\varphi(\psi_1(x_{11}, \ldots, x_{n_1 1}), \ldots, \psi_1(x_{1k}, \ldots, x_{n_1 k})), \ldots, \varphi(\psi_m(x_{N1}, \ldots, x_{M1}), \ldots, \psi_m(x_{Nk}, \ldots, x_{Mk}))\big)$$

$$= \varphi\big(\varphi_i(\psi_1(x_{11}, \ldots, x_{n_1 1}), \ldots, \psi_m(x_{N1}, \ldots, x_{M1})), \ldots, \varphi_i(\psi_1(x_{1k}, \ldots, x_{n_1 k}), \ldots, \psi_m(x_{Nk}, \ldots, x_{Mk}))\big)$$

$= \varphi\big(\psi_a(x_{11}, \ldots, x_{M1}), \ldots, \psi_a(x_{1k}, \ldots, x_{Mk})\big)$. Thus any composition $\psi_a$ commutes with any $\varphi \in \Phi$.

If we proceed now by induction on the number of factors $\varphi_i \in \Phi$ in the composition forming $\psi_b$, a computation identical to the last one shows that $\psi_a$ and $\psi_b$ must commute. ∎

**3.5 Corollary.** — Multiplication in $A$ is commutative and distributes over $\oplus$ from the left.
**Proof:** Both commutativity and distributivity are merely instances of the last Lemma 3.4. Indeed, if we replace $X$ with $A$, then Lemma 3.4 says that *all operations built out of iterations of $\oplus$ commute with one another*. For elements $a, b \in A$, this says that $ab = ba$. For $b(t) = b_1 \oplus b_2 \in A$, this says that

$$a(b_1 \oplus b_2) = b_1 a \oplus b_2 a = ab_1 \oplus ab_2.$$ ∎

**3.6** In summary: we see that $A$ is analogous to $\mathbb{N}$ insofar as both are commutative, multiplicative monoids generated by their multiplicative unit 1 and a binary operation $\oplus$ over which multiplication distributes.

## 4 — Arithmetic in $A$.

We immediately want to know if there are further similarities between $A$ and $\mathbb{N}$. For instance, we'd like to know more about multiplication in $A$, specifically, how factorization behaves in $A$.

To this end, observe that the *magnitude* of any element $a$ in $A$, that is, the number of 1's appearing in any expression in $M$ representing $a$, provides us with a function $|\cdot| : A \to \mathbb{N}$ satisfying

$$|a \oplus b| = |a| + |b|,$$
$$|1| = 1, \quad \text{and} \quad |ab| = |a||b|, \tag{$**$}$$

which is to say that magnitude constitutes a homomorphism from $A$ to $\mathbb{N}$. The multiplicative capacity of this homomorphism gives us one significant factorization property in $A$ right away.

**4.1 Claim.** — Every $a \in A$ admits at least one factorization into irreducibles, meaning that we can write $a$ as a product $a = a_1 \cdots a_n$ in $A$, where each $a_i$ <u>cannot</u> be written as a product in $A$.
**Proof:** Every process of factoring $a$ in $A$ must terminate, since $|\cdot|$ takes any factorization of $a$ in $A$ into a factorization of $|a|$ in $\mathbb{N}$. ∎

**4.2** On the other hand, this factorization need not be unique. In fact, factorization in $A$ is not even cancellative, which is to say that $ab_1 = ab_2$ need not imply that $b_1 = b_2$.

**4.3 Claim.** — Factorization in $A$ is non-cancellative, and consequently non-unique.
**Counter Example:** One example of a non-cancellative product in $A$ is

$$2 \cdot ((1 \oplus 2) \oplus (4_1 \oplus 1)) = 2 \cdot ((1 \oplus 3_1) \oplus (3_2 \oplus 1))$$

(see **(table 1)** above for the definitions of the elements $2, 4_1, 3_1, 3_2 \in A$). Since no expression for any one of the elements $2, 3_1, 3_2, 4_1 \in A$ admits an elementary instance of self-commutativity, the only elementary instances of self-commutativity of $(1 \oplus 2) \oplus (4_1 \oplus 1)$ and $(1 \oplus 3_1) \oplus (3_2 \oplus 1)$ are

$$(1 \oplus 2) \oplus (4_1 \oplus 1) = (1 \oplus 4_1) \oplus (2 \oplus 1) \quad \text{and} \quad (1 \oplus 3_1) \oplus (3_2 \oplus 1) = (1 \oplus 3_2) \oplus (3_1 \oplus 1),$$

respectively. Thus $(1 \oplus 2) \oplus (4_1 \oplus 1) \neq (1 \oplus 3_1) \oplus (3_2 \oplus 1)$. Yet the following sequence of transpositions shows that when we multiply both elements by $2 \in A$, the products become identical:

$$2 \cdot ((1\oplus 2) \oplus (4_1\oplus 1)) =$$
$$= ((1\oplus 2) \oplus (4_1\oplus 1)) \oplus ((1\oplus 2) \oplus (4_1\oplus 1)) =$$
$$= ((1\oplus 2) \oplus (1\oplus 2)) \oplus ((4_1\oplus 1) \oplus (4_1\oplus 1)) =$$
$$= ((1\oplus 1) \oplus (2\oplus 2)) \oplus ((4_1\oplus 1) \oplus (4_1\oplus 1)) =$$
$$= ((1\oplus 1) \oplus (4_1\oplus 1)) \oplus ((2\oplus 2) \oplus (4_1\oplus 1)) =$$
$$= ((1\oplus 1) \oplus (4_1\oplus 1)) \oplus ((2\oplus 4_1) \oplus (2\oplus 1)) =$$
$$= ((1\oplus 1) \oplus (4_1\oplus 1)) \oplus ((3_1\oplus 3_1) \oplus (2\oplus 1)) =$$
$$= ((1\oplus 1) \oplus (3_1\oplus 3_1)) \oplus ((4_1\oplus 1) \oplus (2\oplus 1)) =$$
$$= ((1\oplus 3_1) \oplus (1\oplus 3_1)) \oplus ((4_1\oplus 2) \oplus (1\oplus 1)) =$$
$$= ((1\oplus 3_1) \oplus (1\oplus 3_1)) \oplus ((3_2\oplus 3_2) \oplus (1\oplus 1)) =$$
$$= ((1\oplus 3_1) \oplus (1\oplus 3_1)) \oplus ((3_2\oplus 1) \oplus (3_2\oplus 1)) =$$
$$= ((1\oplus 3_1) \oplus (3_2\oplus 1)) \oplus ((1\oplus 3_1) \oplus (3_2\oplus 1)) =$$
$$= 2 \cdot ((1\oplus 3_1) \oplus (3_2\oplus 1))$$

∎

## 5 — "$q$-Deformations" of Addition.

We would like finally to explain how we located the above instance of non-cancellativity in $A$.

Taking as cue the way in which the morphism $|\cdot| : A \to \mathbb{N}$ gave us information about arithmetic in $A$, we seek an invariant of expressions in $M$ that is preserved under elementary instances of self-commutativity, i.e., one that descends to $A$, yet one that's finer than magnitude.

**5.1 Definition.** — Let $\mathbb{Z}[q]$ be the ring of integral polynomials in the variable $q$. We define a new binary operation $\dotplus : \mathbb{Z}[q] \times \mathbb{Z}[q] \to \mathbb{Z}[q]$ according to

$$f \dotplus g = f + q\, g.$$

This operation $\dotplus$ is self-commutative, and we let $(\mathbb{Z}[q], \dotplus)$ denote the set $\mathbb{Z}[q]$ equipped with $\dotplus$.

**5.2 Claim.** — There exists a unique morphism $(A, \oplus) \longrightarrow (\mathbb{Z}[q], \dotplus)$ taking $A \ni 1 \mapsto 1 \in \mathbb{Z}[q]$, i.e., a unique function $\ell : (A, \oplus) \to (\mathbb{Z}[q], \dotplus)$ such that

$$\ell(1) = 1, \quad \text{and} \quad \ell(a \oplus b) = \ell(a) \dotplus \ell(b).$$

**Proof:** We can take the condition $\ell(a \oplus b) = \ell(a) \dotplus \ell(b) = \ell(a) + q\,\ell(b)$ as the inductive definition of $\ell$. We need only establish that it's well defined. But this is clear since self-commutativity of $\dotplus$ means that every elementary instance of self-commutativity in $M$ is preserved by $\ell$. ∎

**5.3 Claim.** — The function $\ell : (A, \oplus) \to (\mathbb{Z}[q], \dotplus)$ is a homomorphism of multiplicative monoids, in particular, as a "$q$-deformation" of $|\cdot|$, the function $\ell$ retains the property (∗∗) of §4:

$$\ell(ab) = \ell(a)\,\ell(b).$$

**Proof:** We argue by induction on magnitude, it being immediate that $\ell(1 \cdot a) = \ell(1)\ell(a)$. Suppose $\ell$ preserves all products between an arbitrary element $b \in A$ and any element of magnitude $\leq n$. Let $a = a_1 \oplus a_2$ have magnitude $n + 1$. Then by distributivity of multiplication over $\oplus$ from the left in $A$, we have

$$\ell(ab) = \ell(a_1b \oplus a_2b) := \ell(a_1b) + q\,\ell(a_2b)$$
$$= \ell(a_1)\ell(b) + q\,\ell(a_2)\ell(b)$$
$$= (\ell(a_1) \dotplus \ell(a_2))\,\ell(b)$$
$$= \ell(a)\ell(b). \qquad\blacksquare$$

**5.4** Observe that magnitude is nothing but the composition $|\cdot| : A \xrightarrow{\ell} \mathbb{Z}[q] \xrightarrow{q=1} \mathbb{Z}$, whose image lies in $\mathbb{N}$.

The function $\ell$ reflects factorization properties of $A$ in $\mathbb{Z}[q]$. For instance, if $\ell(a)$ is irreducible in $\mathbb{Z}[q]$, then $a$ must be irreducible in $A$.

We can also use $\ell$ to locate candidate elements in $A$ that may have factorization properties that don't hold in $\mathbb{Z}[q]$. For instance, since $\mathbb{Z}[q]$ is a unique factorization domain, if $ab_1 = ab_2$ is a non-cancellative product in $A$, then we must have $\ell(b_1) = \ell(b_2)$ in $\mathbb{Z}[q]$, and elements identified by $\ell$ become candidates for factors in non-cancellative products in $A$.

For this reason, we want $\ell$ to be as fine as possible. But already it's not hard to see that as defined, $\ell$ is not a complete invariant: for instance $\ell(2 \oplus 2) = 1 + 2q + q^2 = \ell(3_1 \oplus 1)$. We can deform $\ell$ a bit more to make it finer. We can make it a complete invariant of all elements of magnitude $\leq 5$ in $A$:

**5.5 Definition.** — Deform $\ell : A \to \mathbb{Z}[q]$ by deforming the operation $\dotplus$ on $\mathbb{Z}[q]$. Specifically, define a new operation $f +_a g := (1 + q)f + (1 - q)g$. This deformed operation $+_a$ remains self-commutative, so continues to define a unique morphism $\ell_a : (A, \oplus) \to (\mathbb{Z}[q], +_a)$, which now becomes

$$\ell_a(a \oplus b) = \ell_a(a) + \ell_a(b) + (\ell_a(a) - \ell_a(b))\,q,$$

and $\ell_a$ retains the properties of Claims 5.2 and 5.3 above.

**5.6** The reader can view this further-deformed function $\ell_a$ as an attempt to measure the non-commutativity of $\oplus$ more accurately by twisting the iterated asymmetries of $\oplus$ antisymmetrically into higher-and-higher powers of $q \in \mathbb{Z}[q]$.

Corresponding to **(table 1)** above, we can fill out a table of the $\ell_a$-values of all elements of magnitude $\leq 5$ in $A$:

| 1 | 2 | 3 | 4 | 5 | |
|---|---|---|---|---|---|
| $\ell_a(1) = 1$ | $\ell_a(2) = 2$ | $\ell_a(3_1) = 3 - q$ | $\ell_a(4_1) = 4$ | $\ell_a(5_1) = 5 - 2q + q^2$ | |
| | | $\ell_a(3_2) = 3 + q$ | $\ell_a(4_2) = 4 - 3q + q^2$ | $\ell_a(5_2) = 5 - q^2$ | |
| | | | $\ell_a(4_3) = 4 + q - q^2$ | $\ell_a(5_3) = 5 + 2q + q^2$ | |
| | | | $\ell_a(4_4) = 4 - q - q^2$ | $\ell_a(5_4) = 5 - 6q + 4q^2 - q^3$ | |
| | | | $\ell_a(4_5) = 4 + 3q + q^2$ | $\ell_a(5_5) = 5 - 2q^2 + q^3$ | |
| | | | | $\ell_a(5_6) = 5 - 2q - 2q^2 + q^3$ | **(table 3)** |
| | | | | $\ell_a(5_7) = 5 + 4q - q^3$ | |
| | | | | $\ell_a(5_8) = 5 - 4q + q^3$ | |
| | | | | $\ell_a(5_9) = 5 + 2q - 2q^2 - q^3$ | |
| | | | | $\ell_a(5_{10}) = 5 - 2q^2 - q^3$ | |
| | | | | $\ell_a(5_{11}) = 5 + 6q + 4q^2 + q^3$ | |
| | | | | $\ell_a(5_{12}) = 5 - 3q$ | |
| | | | | $\ell_a(5_{13}) = 5 + 3q$ | |

**5.7** As described in 5.4 above, we're primarily interested in the question of $\ell_a$'s injectivity. Looking for specific ways in which injectivity of $\ell_a$ can fail, we eventually come upon the following observation.

**5.8 Observation.** — Elements $a_1, a_2, b_1, b_2 \in A$ satisfy
$$\ell_a((1 \oplus a_1) \oplus (a_2 \oplus 1)) = \ell_a((1 \oplus b_1) \oplus (b_2 \oplus 1))$$

if and only if

$$\ell_a(a_1) + \ell_a(a_2) = \ell_a(b_1) + \ell_a(b_2).$$

**Proof:** Expanding $\ell_a((1 \oplus a_1) \oplus (a_2 \oplus 1))$, we find that

$$\ell_a((1 \oplus a_1) \oplus (a_2 \oplus 1)) = 2 + \ell_a(a_1) + \ell_a(a_2) - (2 + \ell_a(a_1) + \ell_a(a_2)) q^2,$$

making the observation immediate. ■

**5.9** A glance at **(table 3)** above reveals that $\ell_a(2) + \ell_a(4_1) = \ell_a(3_1) + \ell_a(3_2)$ is an instance of Lemma 5.3. Hence

$$\ell_a((1 \oplus 2) \oplus (4_1 \oplus 1)) = \ell_a((1 \oplus 3_1) \oplus (3_2 \oplus 1)),$$

and the distinct elements $(1 \oplus 2) \oplus (4_1 \oplus 1)$ and $(1 \oplus 3_1) \oplus (3_2 \oplus 1)$ in $A$ become candidates for a possible non-cancellative product in $A$. The first multiplication we try gives us the counter example

$$2 \cdot ((1 \oplus 2) \oplus (4_1 \oplus 1)) = 2 \cdot ((1 \oplus 3_1) \oplus (3_2 \oplus 1)).$$